\documentclass[final,1p,times]{elsarticle}
\usepackage{amsmath,amsfonts,amssymb,amsthm}
\usepackage[all]{xy}
\usepackage{graphicx}
\usepackage{aliascnt,hyperref}

\journal{arXiv}

\newtheorem{teo}{Theorem}

\def\partet#1#2#3#4{\newaliascnt{#1}{#2}\newtheorem{#1}[#1]{#3}\aliascntresetthe{#1}\providecommand*{#4}{#3}}
\def\parted#1#2#3#4{\newaliascnt{#1}{#2}\newdefinition{#1}[#1]{#3}\aliascntresetthe{#1}\providecommand*{#4}{#3}}
\partet{lema}{teo}{Lemma}{\lemaautorefname}
\partet{prop}{teo}{Proposition}{\propautorefname}
\partet{cor}{teo}{Corollary}{\corautorefname}
\parted{defs}{teo}{Definition}{\defsautorefname}
\parted{ejemplo}{teo}{Example}{\ejemploautorefname}
\parted{obs}{teo}{Remark}{\obsautorefname}
\parted{parr}{teo}{}{\parrautorefname}
\parted{notc}{teo}{Notation}{\notcautorefname}
\newproof{dem}{Proof}
\newtheorem*{nteo}{Theorem}

\def\Z{\mathbb{Z}}
\def\C{\mathbb{C}}
\def\N{\mathbb{N}}

\def\PP{\mathbb{P}}
\def\R{\mathbb{R}}

\def\then{\Longrightarrow}

\def\fl#1{\stackrel{#1}{\longrightarrow}}

\begin{document}

\begin{frontmatter}

\title{Algorithm to find a maximum of a multilinear map over a product of spheres.}
\author[dm]{C\'esar Massri\corref{correspondencia}\fnref{finanaciado}}
\address[dm]{Department of Mathematics, FCEN, University of Buenos Aires, Argentina}
\cortext[correspondencia]{Address for correspondence: Department of Mathematics, FCEN, University of Buenos Aires, Argentina}
\fntext[finanaciado]{The author was fully supported by CONICET, Argentina}
\ead{cmassri@dm.uba.ar}

\begin{abstract}
We provide an algorithm to compute the 2-norm maximum of a multilinear map
over a product of spheres.
As a corollary we give a method to compute the first singular value of a linear map and an
application to the theory of entangled states in quantum physics. Also,
we give an application to find the closest rank-one tensor of a given one.
\end{abstract}

\end{frontmatter}

\section*{Introduction.}
A lot of problems in mathematics need to maximize a bilinear form over a product of spheres,
for example the 2-norm of a matrix is given by the maximum of the bilinear form $(x,y)\rightarrow x^tAy$,
where $\|x\|=\|y\|=1$. Another interesting problem is to find the closest rank-one tensor
of a given tensor $\sum a_{ijk}x_i\otimes y_j\otimes z_k$. To answer this problem one has to
find the maximum of a trilinear form over a product of three spheres (see the examples).\\

This article provides an algorithm to find the maximum of a multilinear map over
a product of spheres,
$$\ell:\R^{n_1+1}\times\ldots\times\R^{n_r+1}\rightarrow\R^{n_{r+1}+1},\quad \max_{\|x_1\|=\ldots=\|x_r\|=1}\|\ell(x_1,\ldots,x_r)\|.$$
We have reduced the problem of finding the maximum of $\ell$ to a problem of finding
fixed points of a map $\nabla\ell:\PP^{n_1}\times\ldots\times\PP^{n_{r+1}}\rightarrow
\PP^{n_1}\times\ldots\times\PP^{n_{r+1}}$.
The advantage of this reduction is the possibility to count the number of extreme points of $\ell$,
and also, to find the fixed points of $\nabla\ell$ solving a system of polynomial equations.
There are standard algebro-geometric tools to solve systems of polynomial equations.\\

In {\bf \autoref{sec-def}} we review some concepts and definitions in algebraic geometry,
such as, projective space, maps, products of projective spaces and maps between them.
We use these definitions in the article.

In {\bf \autoref{sec-theo}}, using Lagrange's method of multipliers, see \cite[\S 13.7]{MR0344384},
we reduce the problem of finding the maximum of
a multilinear map $\ell$, to the problem of finding fixed points of a map $\nabla\ell$.
We compare our approach with the ones in the literature.

In {\bf \autoref{sec-num}} we make a digression to discuss the number of extreme points of a multilinear map
over a product of spheres. We use intersection theory to count the number of fixed points of the
map $\nabla\ell:\PP^{n_1}\times\ldots\times\PP^{n_{r+1}}\rightarrow
\PP^{n_1}\times\ldots\times\PP^{n_{r+1}}$. Recall that the number
of fixed points of a generic map $F:\PP^N\rightarrow\PP^N$ of degree $d$ is
$1+d+\ldots+d^N$. In this section we give a formula to compute the number of extreme points of a multilinear map
over a product of spheres. If the map is generic, this number is achieved over $\C$, and if it is not generic, this
number is a bound when the extreme points are finite.
In the literature, the extreme points of $\ell$ are called singular vectors (see \cite{1574201}) and
in this section we count them.

In {\bf \autoref{sec-bil}} we use our approach to find the maximum
of a bilinear form over a product of spheres. In the bilinear case, the map $\nabla\ell$,
induces a linear map $L:\PP^N\rightarrow\PP^N$, where $N$ is a natural number,
and we prove that for a generic $q\in\PP^N$,
the sequence $\{q,L(q),L^2(q),\ldots\}$ converges to the absolute maximum. In other words,
the absolute maximum is an attractive fixed point of $L$.
Also, with the same tools, we give an algorithm to find the spectral radius of a square matrix.

In {\bf \autoref{sec-pres}} we use the theory developed to present the algorithm. We take
advantage of a result in \autoref{sec-theo}; the classes of extreme points of the multilinear form $\ell$ are in
bijection with the fixed points of the map $\nabla\ell$.
We reduce the problem of finding fixed points of $\nabla\ell$ to
solve a system of polynomial equations with finitely many solutions.
In the literature about computational aspects of algebraic geometry, 
there exists a lot of algorithms to solve a system of polynomial equations with finitely many solutions, see \cite{MR2161984}.
This gives us the ability to find the absolute maximum of $\ell$.
It is important to mention that the system of polynomial equations obtained with our approach is
slightly different from  the system of polynomial
equations obtained naively from the method of Lagrange's multipliers.
Our approach in projective geometry, allows us to find the correct solution removing some constrains.
In the first part of the section, we present a direct method to find the maximum value
of a generic multilinear form over a product of spheres. Basically, it reduces to find the spectral radius of a matrix.
In the second part of the section, we give an algorithm to find the 
point $(x_1,\ldots,x_r)\in\R^{n_1+1}\times\ldots\times\R^{n_r+1}$, 
where $\|x_i\|=1$, $1\leq i\leq r$, such that $|\ell(x_1,\ldots,x_r)|$ is maximum.
This last algorithm, requires the extra hypothesis, $2n_1,\ldots,2n_r\leq n_1+\ldots+n_r$.

In {\bf \autoref{sec-ex}} we use the theory developed to compute a lot of examples and applications.
One of them is the ability to find the closest rank-one tensor of a given tensor.
We prove that this problem is well posed and we apply our algorithm to solve it.
Another application that we will give is related to quantum physics.
It is a criterion of separability, given a quantum state,
we can say if it is separable (see \autoref{ex-quantum} for definitions and related concepts).

\section{Review on Projective Geometry.}\label{sec-def}
In this section we give basics definitions that we are going to use,
such as, projective space, maps, projective tangent space, product of projective spaces and maps between them.
In this section we are assuming
that the base field is $\R$, but all the definitions are true in the complex case.
All the notions in this section may be found in \cite{MR1182558}.

\begin{defs}
Let $n$ be a natural number and let $\R^{n+1}$ be a real vector space of dimension $n+1$.
The \emph{projective space}, $\PP^n$, is the space of lines passing throw the origin in $\R^{n+1}$.
We say that the \emph{dimension} of $\PP^n$ is $n$.
Every nonzero vector $v$ in $\R^{n+1}$ determines the line $[v]$ that joins $v$ with the origin $0\in\R^{n+1}$.
The vector $v$ and $\lambda v$, for $\lambda\in\R$, $\lambda\neq 0$, determines the same point $[v]\in\PP^n$.\\

Let's fix a basis $\{v_0,\ldots,v_n\}$ of $\R^{n+1}$.
If the coordinates, in this basis, of $v$ are $(a_0,\ldots,a_n)$, then the coordinates of the point $[v]$ are
$$[v]=(a_0:\ldots:a_n)=(\lambda a_0:\ldots:\lambda a_n),\quad\lambda\in\R,\,\lambda\neq 0.$$
In general we denote $[v]\in\PP^n$ to remark that the point $[v]$ is represented by the
vector $v\in\R^{n+1}$. Also, we denote an arbitrary point in projective space, $p\in\PP^n$.
The projective space $\PP^n$ is a compact space.\\

Let $n$ and $m$ be two natural numbers. We say that a polynomial $P$ in $n+1$ variables
is \emph{homogeneous} of \emph{degree} $d$, where $d$ is a natural number, if
$$P(\lambda x_0,\ldots,\lambda x_n)=\lambda^d P(x_0,\ldots,x_n),\quad \lambda\in\R,\,\lambda\neq 0.$$
For example, a linear form is homogenous of degree $1$.\\

A map $F$ from $\PP^n$ to $\PP^m$, denoted $F:\PP^n\rightarrow\PP^m$, is given by $m+1$ homogeneous
polynomials, $F_0,\ldots,F_m$, of degree $d$
$$F=(F_0:\ldots:F_m):\PP^n\rightarrow\PP^m,\quad F(x)=(F_0(x):\ldots:F_m(x)),\,x\in\PP^n.$$
The homogeneity of the polynomials $F_0,\ldots,F_m$, implies that
the value of $F$ at $[v]$ and at $[\lambda v]$ is the same in $\PP^m$.
We say that $F$ has \emph{degree} equal to $d$. When $d=1$ we say that $F$ is \emph{linear}.\\

Let $n_1,\ldots,n_{r}$ be a list of natural numbers.
A \emph{multihomogeneous polynomial} is a polynomial $P$ in variables $x_0^i,\ldots,x_{n_i}^i$,
for $1\leq i\leq r$, such that
$$P(\lambda_1x^1,\ldots,\lambda_rx^r)=\lambda_1^{d_1}\ldots\lambda_r^{d_r}P(x^1,\ldots,x^r),\quad x^i=(x_0^i,\ldots,x_{n_i}^i).$$
The vector $(d_1,\ldots,d_r)$ is called the \emph{multidegree} of $P$.
For example, a multilinear form is a multihomogeneous polynomial of multidegree $(1,\ldots,1)$.\\

A map $F:\PP^{n_1}\times\ldots\times\PP^{n_r}\rightarrow \PP^{m}$, where $m\in\N$,
is given by $m+1$ multihomogeneous polynomials, $F_0,\ldots,F_m$, of multidegree $(d_1,\ldots,d_r)$,
$$F(x^1,\ldots,x^r)=(F_0(x^1,\ldots,x^r):\ldots:F_m(x^1,\ldots,x^r)),\quad x^i\in\PP^{n_i},\,1\leq i\leq r.$$
The multi-homogeneity of the polynomials $F_0,\ldots,F_m$, implies that
the value of $F$ at $([v^1],\ldots,[v^r])$ and at $([\lambda_1v^1],\ldots,[\lambda_rv^r])$ is the same in $\PP^m$.
We say that $F$ has \emph{multidegree} $(d_1,\ldots,d_r)$.\\

Finally, a map $F:\PP^{n_1}\times\ldots\times\PP^{n_r}\rightarrow \PP^{m_1}\times\ldots\times\PP^{m_s}$
is given by $s$ maps $F=(F_1,\ldots,F_s)$,
$$F_i:\PP^{n_1}\times\ldots\times\PP^{n_r}\rightarrow\PP^{m_i},\quad 1\leq i\leq s.$$
Note that the multidegree of $F_i$ may differs from the multidegree of $F_j$, $i\neq j$.
When all the forms $\{F_1,\ldots,F_s\}$ are multilinear, we say that $F$ is a \emph{multilinear} map.
\end{defs}

\begin{defs}
Let $n$ and $m$ be two natural numbers and fix a basis for $\R^{n+1}$ and for $\R^{m+1}$.
Every vector $v\in\R^{n+1}$ has associated a vector space of dimension $n+1$;
the \emph{tangent space}, denoted $T_v\R^{n+1}$.

A polynomial map $F=(F_0,\ldots,F_m):\R^{n+1}\rightarrow\R^{m+1}$ such that $F(v)=w$ determines a linear map, $\widehat{dF_v}$,
called the differential of $F$ at $v$,
$$\widehat{dF_v}:T_v\R^{n+1}\rightarrow T_w\R^{m+1},\quad
\widehat{dF_v}(a_0,\ldots,a_n)=\left(\sum_{i=0}^n\frac{\partial F_0}{\partial x_i}(v)a_i,\ldots,
\sum_{i=0}^n\frac{\partial F_m}{\partial x_i}(v)a_i\right).$$

In projective space the situation is similar, \cite[p.181]{MR1182558}.
Every point $x\in\PP^n$ has associated an $n$-dimensional
projective space; the \emph{projective tangent space}, denoted $\mathbb{T}_x\PP^n$.
A map $F=(F_0:\ldots:F_m):\PP^n\rightarrow\PP^m$ of degree $d$ such that $F(x)=y$
induces a linear map between projective tangent spaces,
$$dF_x:\mathbb{T}_x\PP^n\rightarrow \mathbb{T}_y\PP^m,\quad
dF_x(a_0:\ldots:a_n)=\left(\sum_{i=0}^n\frac{\partial F_0}{\partial x_i}(x)a_i:\ldots:
\sum_{i=0}^n\frac{\partial F_m}{\partial x_i}(x)a_i\right).$$
Given that the partial derivative of a homogeneous polynomial is also homogeneous,
the map $dF_x$ is well defined.
\end{defs}

\begin{obs}\label{def-euler}
Recall the \emph{Euler relation} for a homogeneous polynomial $P$ of degree $d$, \cite[p. 182]{MR1182558},
$$\sum_{i=0}^N \frac{\partial P}{\partial x_i}(v)v_i=d\cdot P(v),\quad v=(v_0,\ldots,v_N)\in\R^{N+1}.$$
The relation follows at once by differentiating both sides of the equation $P(\lambda v)=\lambda^d P(v)$.

If $F=(F_0:\ldots:F_N):\PP^N\rightarrow\PP^N$ is a map of degree $d$ and $x\in\PP^N$ is a point such that $F(x)=x$,
then, using the Euler relation, we get $dF_x(x)=x$,
$$dF_x:\mathbb{T}_x\PP^N\rightarrow \mathbb{T}_x\PP^N,\quad
dF_x(a_0:\ldots:a_N)=\left(\sum_{i=0}^N\frac{\partial F_0}{\partial x_i}(x)a_i:\ldots:
\sum_{i=0}^N\frac{\partial F_N}{\partial x_i}(x)a_i\right).$$
In particular, if the vector $v\in\R^{N+1}$ represents $x\in\PP^N$, $x=[v]$,
and the matrix $\widehat{dF_x}$ represent the linear map $dF_x$,
$$\left(\widehat{dF_x}\right)_{i+1,j+1}=\frac{\partial F_i}{\partial x_j}(v),\quad 0\leq i,j\leq N,$$
then, $v$ is an eigenvector of $\widehat{dF_x}$.
Let's compute the eigenvalue of the eigenvector $v$.
Given that $F(x)=x$ there exists a nonzero real number $\lambda$ such that
$(F_0(v),\ldots,F_N(v))=\lambda v$. Then
$$\lambda v_j=F_j(v)=\frac{1}{d}\sum_{i=0}^N \frac{\partial F_j}{\partial x_i}(v)v_i,\quad 0\leq j\leq N.$$
Then, the eigenvalue of $v$ is $d\cdot\lambda$, where $d$ is the degree of the map $F$.
\end{obs}

\section{Theory for a multilinear map.}\label{sec-theo}

In this section we translate the problem of finding a maximum of a multilinear map
to a problem of finding fixed points. Let's present the notation and some basics preliminaries.

Let $\mathbb{S}^{n}$ be the sphere in $\R^{n+1}$,
$$\mathbb{S}^{n}=\{u\in\R^{n+1}\colon\,\|u\|=\sqrt{|u_0|^2+\ldots+|u_n|^2}=1\},$$
and let $\langle-,-\rangle:\R^{n+1}\times\R^{n+1}\rightarrow\R$ be the inner product, $\langle x,y\rangle=x_0y_0+\ldots+x_ny_n$.
The norm associated to this inner product is the usual 2-norm, $\langle u,u\rangle=\|u\|^2$.

When the codomain of a map is $\R$, we say that the map is a \emph{form}.
\begin{lema}
Given a multilinear map $\ell:\R^{n_1+1}\times\ldots\times\R^{n_r+1}\fl{}\R^{s+1}$ there exists a multilinear form $\hat{\ell}$,
$$\hat{\ell}:\R^{n_1+1}\times\ldots\times\R^{n_r+1}\times\R^{s+1}\fl{}\R,\quad
\hat{\ell}(x_1,\ldots,x_r,y)=\langle\ell(x_1,\ldots,x_r),y\rangle,$$
such that
$$\max_{\|x_1\|=\ldots=\|x_r\|=1}\|\ell(x_1,\ldots,x_r)\|=\max_{\|x_1\|=\ldots=\|x_r\|=\|y\|=1}|\hat{\ell}(x_1,\ldots,x_r,y)|.$$
\end{lema}
\begin{dem}
The proof is bases on the compactness of the sphere.
Let $(x_1,\ldots,x_r)\in\mathbb{S}^{n_1}\times\ldots\times\mathbb{S}^{n_r}$ be a point such that
$z=\ell(x_1,\ldots,x_r)$ has the maximum norm and let $y=z/\|z\|$.
Then
$$|\hat{\ell}(x_1,\ldots,x_r,y)|=|\langle z,y\rangle|=\frac{\langle z,z\rangle}{\|z\|}=
\|z\|=\|\ell(x_1,\ldots,x_r)\|\then$$
$$\max_{\|x_1\|=\ldots=\|x_r\|=\|y\|=1}|\hat{\ell}(x_1,\ldots,x_r,y)|\geq
\max_{\|x_1\|=\ldots=\|x_r\|=1}\|\ell(x_1,\ldots,x_r)\|.$$
Analogously, let $(x_1,\ldots,x_r,y)\in\mathbb{S}^{n_1}\times\ldots\times\mathbb{S}^{n_r}\times\mathbb{S}^{s}$ be
a point such that $|\hat{\ell}(x_1,\ldots,x_r,y)|$ is maximum.
Let $z=\ell(x_1,\ldots,x_r)$. Then,
$$|\hat{\ell}(x_1,\ldots,x_r,y)|=|\langle z,y\rangle|\leq\|z\|\|y\|=\|\ell(x_1,\ldots,x_r)\|\then$$
$$\max_{\|x_1\|=\ldots=\|x_r\|=\|y\|=1}|\hat{\ell}(x_1,\ldots,x_r,y)|\leq
\max_{\|x_1\|=\ldots=\|x_r\|=1}\|\ell(x_1,\ldots,x_r)\|.$$\qed
\end{dem}

As a corollary of the previous lemma, we will work with multilinear forms. Specifically,
to make the notation easiest, we will work with $\ell:\R^{n+1}\times\R^{m+1}\times\R^{s+1}\fl{}\R$ a trilinear form.
Our goal is to find the maximum of $\ell$ over a product of three spheres.

Using Lagrange's method of multipliers, (\cite[\S 13.7]{MR0344384}), we know that the
extreme points of $\ell$, over $\mathbb{S}^{n}\times\mathbb{S}^{m}\times\mathbb{S}^{s}$, satisfy
$$\left\{\begin{array}{lcr}
\partial\ell/\partial x_i(x_0,\ldots,x_n,y_0,\ldots,y_m,z_0,\ldots,z_s)&=&2\alpha x_i,\quad 0\leq i\leq n,\\
\partial\ell/\partial y_j(x_0,\ldots,x_n,y_0,\ldots,y_m,z_0,\ldots,z_s)&=&2\beta y_j,\quad 0\leq j\leq m,\\
\partial\ell/\partial z_k(x_0,\ldots,x_n,y_0,\ldots,y_m,z_0,\ldots,z_s)&=&2\lambda z_k,\quad 0\leq k\leq s,
\end{array}\right.$$
$$\quad\alpha,\beta,\lambda\in\R,\quad\|x\|=\|y\|=\|z\|=1.$$
Let's use a better notation,
$$x=(x_0,\ldots,x_n),\,
y=(y_0,\ldots,y_m),\,z=(z_0,\ldots,z_s),$$
$$\frac{\partial\ell}{\partial x}(x,y,z)=
\left(\frac{\partial\ell}{\partial x_0}(x,y,z),\ldots,\frac{\partial\ell}{\partial x_n}(x,y,z)\right),\quad
\frac{\partial\ell}{\partial y}(x,y,z)=\left(\frac{\partial\ell}{\partial y_0}(x,y,z),\ldots,\frac{\partial\ell}{\partial y_m}(x,y,z)\right),$$
$$\frac{\partial\ell}{\partial z}(x,y,z)=\left(\frac{\partial\ell}{\partial z_0}(x,y,z),\ldots,
\frac{\partial\ell}{\partial z_s}(x,y,z)\right),\quad
\nabla\ell(x,y,z)=\left(\frac{\partial\ell}{\partial x}(x,y,z),\frac{\partial\ell}{\partial y}(x,y,z),
\frac{\partial\ell}{\partial z}(x,y,z)\right).$$

\begin{defs}
A point $(x,y,z)\in\mathbb{S}^{n}\times\mathbb{S}^{m}\times\mathbb{S}^{s}$ is called
an \emph{extreme point} of $\ell$ if it satisfies the system of equations
$$\nabla\ell(x,y,z)=(2\alpha x,2\beta y,2\lambda z),$$
for some $\alpha,\beta,\lambda\in\R$.
Note that if $(x,y,z)$ is an extreme point, then $(\pm x,\pm y,\pm z)$ is also an extreme point.
We say that they belong to the same \emph{class}.
\end{defs}

\begin{prop}\label{theo-bij}
There is a bijection between classes of extreme points of $\ell$ and fixed points of the map
$$\nabla\ell:
\PP^n\times\PP^m\times\PP^s\rightarrow\PP^n\times\PP^m\times\PP^s,$$
$$([x],[y],[z])\rightarrow
\left(\frac{\partial\ell}{\partial x}([x],[y],[z]),
\frac{\partial\ell}{\partial y}([x],[y],[z]),\frac{\partial\ell}{\partial z}([x],[y],[z])\right).$$
\end{prop}
\begin{dem}
Given an extreme point $(x,y,z)$, consider $([x],[y],[z])$. This assignment is independent of the class of $(x,y,z)$.
By definition, it gives a fixed point of $\nabla\ell$.

Given a fixed point $([x],[y],[z])\in\PP^n\times\PP^m\times\PP^s$ of $\nabla\ell$, consider representatives $x,y,z$
such that $\|x\|=\|y\|=\|z\|=1$. Then $(x,y,z)$ is an extreme point of $\ell$.\qed
\end{dem}

\begin{obs}\label{theo-hyp}
The map $\nabla\ell:\PP^n\times\PP^m\times\PP^s\rightarrow\PP^n\times\PP^m\times\PP^s$
from \autoref{theo-bij} is not defined over the closed subset
$$\left\{(x,y,z)\,|\,\frac{\partial\ell}{\partial x}(x,y,z)=0\text{ or }
\frac{\partial\ell}{\partial y}(x,y,z)=0\text{ or }
\frac{\partial\ell}{\partial z}(x,y,z)=0
\right\}\subseteq \PP^n\times\PP^m\times\PP^s,$$
but this set is empty if and only if the \emph{hyperdeterminant} of $\ell$ is zero.
The hyperdeterminant is a polynomial in the coefficient of $\ell$,
for the definition and some properties see \cite[\S14]{MR1264417}.

By a result in \cite[\S14, 1.3]{MR1264417}, if
$$2n,2m,2s\leq n+m+s$$
then a generic choice of $\ell$ will make $\nabla\ell$ defined everywhere.
\end{obs}

In the article \cite{1574201}, there is a definition of singular values and singular vectors
for a multilinear form $\ell$. For example, for a trilinear form $\ell$,
the author defined the singular vectors of $\ell$ as the solutions of the
system $\nabla\ell(x,y,z)=(2\alpha x,2\beta y,2\lambda z)$. It is the same as our definition of extreme points.
It is of interest to know the number of singular values/vectors of $\ell$, and in \autoref{sec-num}, we count them.
In the same article, the author proved that the first
singular value is the maximum of $\ell$ over a product of spheres. Also,
under the hypothesis $2n,2m,2s\leq n+m+s$, he proved that
the \emph{hyperdeterminant} of $\ell$ is zero if and only if $0$ is a singular value of $\ell$.
The hyperdeterminant is a polynomial in the coefficients of $\ell$, so, in the generic case,
the number $0$ is not a singular value of $\ell$.\\
There exists another article to mention, \cite{MR1780272}. In it, the authors proposed a multidimensional
singular value decomposition, but it does not preserve the properties that we need, for example, that
the first singular value of $\ell$ corresponds to the maximum of $\ell$ over a product of spheres.

\section{Number of extreme points of a multilinear form.}\label{sec-num}
In this section we use Intersection Theory (\cite[8.4]{MR1644323}) to count the number of fixed
points of a generic map
$\PP^n\times\PP^m\times\PP^s\rightarrow\PP^n\times\PP^m\times\PP^s$ over $\C$.
Recall from \autoref{theo-bij} that there is a bijection between fixed points of
$$\nabla\ell:\PP^n\times\PP^m\times\PP^s\rightarrow\PP^n\times\PP^m\times\PP^s$$
and classes of extreme points of the trilinear form $\ell$
over $\mathbb{S}^n\times\mathbb{S}^m\times\mathbb{S}^s$.
The number of fixed points of $\nabla\ell$
gives a bound to the number of classes of extreme points of $\ell$ that
contains a point with maximum value.
It is known that if $F:\PP^{N}\rightarrow\PP^{N}$ is a generic map of degree $d$, then $F$
has $1+d^2+\ldots+d^{N}$ fixed points, \cite[1.3]{MR2648690}.
Here we generalize this result to a generic map between products of projective spaces.\\

Before we continue with this section, let's make a survey of some related concepts that are in the literature.

In \cite{MR1931400}, \cite{MR2270090}, \cite{MR2296920}, \cite{Cartwright2011}, \cite{MR2854605} and \cite{BKP11}
there is a notion of eigenvectors and eigenvalues associated to a multilinear form $\ell$.
There are a lot of applications and in \cite{Cartwright2011}, the authors counted the number of
eigenvalues of $\ell$ as the number of roots of a characteristic polynomial associated to $\ell$. The idea is
to look at $\ell:\C^n\times\ldots\times\C^n\rightarrow\C$ as a polynomial map $P:\C^n\rightarrow(\C^n)^\vee\cong\C^n$,
$P(x)=\ell(x,\ldots,x,-)$ and then, an eigenvector of $\ell$ is a vector $x\in\C^n$ such that $P(x)=\lambda x$.
If $\ell$ is $m$-multilinear, $P$ has degree $m-1$ and as a map $\PP^{n-1}\rightarrow\PP^{n-1}$
it has $(m-1)^{n-1}+(m-1)^{n-2}+\ldots+1$ fixed points, i.e. eigenvectors of $P$. They arrived at
this number using toric varieties and Newton polytopes.\\

In \cite[7.1.4]{MR2161984} and \cite[3.1]{MR1342295} there is a theory
of multihomogeneous B\'ezout number, or $m$-B\'ezout.
The $m$-B\'ezout gives an upper bound on the cardinality of the intersection
of multihomogeneous polynomials in $\PP^{n_1}\times\ldots\times\PP^{n_k}$.
Given that we are counting the fixed points of a
map $F:\PP^{n_1}\times\ldots\times\PP^{n_k}\rightarrow \PP^{n_1}\times\ldots\times\PP^{n_k}$,
in order to apply this formula, we need to realize the fixed points of $F$ as
an intersection in some product of projective spaces. Concretely,
the intersection of the graph of $F$ and the diagonal.
Let's make an explicit example. Assume for simplicity, that $F$ is linear, $F:\PP^n\rightarrow \PP^n$, we will see that
the $m$-B\'ezout formula gives a very bad bound. Recall that the number of fixed points in this
situation is the number of eigenvectors, that is, $n+1$.
Let's apply the formula to the equations of the graph $\Gamma=\{(x,F(x))\}$ and the diagonal $\triangle=\{(x,x)\}$.
The points in the intersection satisfy the following equations,
$$((x_0:\ldots:x_n),(y_0:\ldots:y_n))\in\PP^n\times\PP^n,\quad
y_iF_j(x)=y_jF_i(x),\quad x_iy_j=x_jy_i,\quad 0\leq i,j\leq n.$$
Note that the equations correspond to the fact that the following matrices have rank one,
$$\begin{pmatrix}
y_0&\ldots&y_n\\
F_0(x)&\ldots&F_n(x)
\end{pmatrix},\quad
\begin{pmatrix}
x_0&\ldots&x_n\\
y_0&\ldots&y_n\end{pmatrix}.$$
By abuse of notation, we denote the equations,
$$(x,y)\in\PP^n\times\PP^n,\quad y=F(x),\quad x=y.$$
Given that the equations have bidegree $(1,1)$,
the $m$-B\'ezout number is the coefficient of $\alpha_1^{n+1}\alpha_2^{n+1}$
in the polynomial $(\alpha_1+\alpha_2)^{2n+2}$. It is the binomial $\binom{2n+2}{n+1}\neq n+1$.

Bernstein proved in \cite{MR0435072} that the number of solutions of a \emph{sparse system} equals
the \emph{mixed volume} of the corresponding \emph{Newton polytopes}. A sparse system is a collection of
Laurent polynomials,
$$f_i=\sum_{(v_1,\ldots,v_n)\in\mathcal{A}_i}c_{i,v_1,\ldots ,v_n}x_1^{v_1}\ldots x_n^{v_n},\quad i=1,\ldots,n$$
where $\mathcal{A}_i$ are fixed finite subsets of $\Z^n$. Its convex hull $Q_i=\text{conv}(\mathcal{A}_i)\subseteq\R^n$ is
called the Newton polytope of $f_i$. Consider the function
$$R(\lambda_1,\ldots,\lambda_n):=\text{vol}(\lambda_1 Q_1+\ldots+\lambda_n Q_n),\quad \lambda_i\geq 0,\,i=1,\ldots,n.$$
where \emph{vol} is the usual Euclidean volume in $\R^n$ and $Q+Q'$ denotes the Minkowski sum of polytopes.
It is a fact that $R$ is a homogeneous polynomial and the coefficient of the monomial $\lambda_1\ldots\lambda_n$
is called
the mixed volume of $Q_1,\ldots,Q_n$. The mixed volume (i.e the number of solutions of a sparse system)
is a very difficult number to compute by hand. In some situations, this is possible and in the general case,
there are a lot of algorithms to compute it. In our situation, we are working with a multihomogeneous polynomial
system, and using Bernstein's theorem, in the paper \cite{MR1720109}, the author gives a recursive formula to compute this number.
In fact, it is proved that, under some hypothesis, if the system is over $\R$ and the functions are generic,
then all the solutions are reals. Here, we present a different and more direct method using intersection theory.\\

In the article \cite[3]{1574201} there is a definition of generalized singular values
for a generic multilinear form. In this section we count them.\\

Let's make an introduction to Intersection Theory.
The germ of intersection theory is the Fundamental Theorem of Algebra. It implies that
given a generic homogeneous polynomial in two variables $F$ of degree $d$, the set of zeroes
$\{x\in\PP^1,\,F(x)=0\}$ has $d$ points. Generalizing this idea, B\'ezout's theorem,
says that given two generic homogeneous polynomials in three variables of degree $d$ and $e$, the set of
zeroes $\{x\in\PP^2,\,F_1(x)=F_2(x)=0\}$ consists of $de$ points. In $\PP^r$ the situation is similar,
if $F_1,\ldots,F_r$ are generic homogeneous polynomials of degree $d_1,\ldots,d_r$ respectively,
the intersection has $d_1d_2\ldots d_r$ points.\\
To formalize these ideas, let's introduce the Chow ring of $\PP^r$, \cite[proof of Prop. 8.4]{MR1644323}
$$A(\PP^r)=\Z[\alpha]/(\alpha^{r+1}).$$
Every variety $X\subseteq \PP^r$ has a class, $[X]\in A(\PP^r)$. The intersection of two
generic varieties $X\cap Y$ corresponds to the product of the classes $[X].[Y]=[X\cap Y]$.
Two different varieties may correspond to the same class, for example, every hypersurface of degree $d$
corresponds to the same class, $d\alpha$, where $\alpha$ is the class of a hyperplane.
For example, $\alpha^r$ corresponds to the intersection of $r$ generic hyperplanes, i.e. a point. The product
$$(d_1\alpha).(d_2\alpha).\ldots.(d_r\alpha)=d_1\ldots d_r \alpha^{r}$$
corresponds to the intersection of $r$ generic hypersurface of degree $d_1,\ldots,d_r$ respectively. We
get $d_1\ldots d_r$ points in the intersection as mentioned.
The class of a variety of codimension $c$ is a homogeneous polynomial of degree $c$ in $\Z[\alpha]/(\alpha^{r+1})$.

The Chow ring is very useful to solve problems in enumerative geometry.
For example, to count the number of fixed points of a generic map $\PP^r\rightarrow\PP^r$
the procedure is the following. Let $A(\PP^r\times \PP^r)$
be the Chow ring of $\PP^r\times \PP^r$, it is $A(\PP^r\times \PP^r)=\Z[a,\alpha]/(a^{r+1},\alpha^{r+1})$, \cite[Ex. 8.4.2]{MR1644323}.
Let $[\triangle]\in A(\PP^r\times \PP^r)$ be
the class of the diagonal, $\triangle=\{(x,x)\}$, and let $[\Gamma]\in A(\PP^r\times \PP^r)$ be
the class of the graph of $F$, $\Gamma=\{(x,F(x))\}$. Given that
$$\dim\triangle+\dim \Gamma=\dim(\PP^r\times \PP^r),$$
the product $[\triangle].[\Gamma]$ is a multiple of the class of a point, $da^r\alpha^r$, \cite[\S8.3]{MR1644323}.
The coefficient $d$ is the number of fixed points of $F$.

The Chow ring of a product of projective spaces, \cite[Ex. 8.3.7]{MR1644323}, is
$$A(\PP^{n_1}\times \ldots\times \PP^{n_k})=A(\PP^{n_1})\otimes_\Z \ldots\otimes_\Z A(\PP^{n_k})=
\Z[\alpha_1,\ldots,\alpha_k]/(\alpha_1^{n_1+1},\ldots,\alpha_k^{n_k+1}).$$
Note that in $A(\PP^{n_1}\times \ldots\times \PP^{n_k})$ there is only one class of a point, $\alpha_1^{n_1}\ldots\alpha_{k}^{n_k}$,
so there is a well defined map called \emph{degree}. The degree of a class is the coefficient
of $\alpha_1^{n_1}\ldots\alpha_{k}^{n_k}$. It may be zero (or negative).\\
The last thing to mention is that every map
$F:\PP^{n_1}\times \ldots\times \PP^{n_k}\rightarrow \PP^{m_1}\times \ldots\times \PP^{m_l}$
induces a morphism of rings, \cite[Prop. 8.3 (a)]{MR1644323},
$$F^\star:A(\PP^{m_1}\times \ldots\times \PP^{m_l})\rightarrow A(\PP^{n_1}\times \ldots\times \PP^{n_k}),\quad
F^\star([X])=[F^{-1}(X)].$$
For a more extensive treatment of intersection theory, see \cite[\S A]{MR0463157}, \cite{MR1644323}.\\

Let's use the previous introduction.
First, we will compute the number of fixed points of a generic map $F:\PP^r\rightarrow\PP^r$ of degree $d$.
Then, we will adapt the proof to a generic map $F:\PP^n\times\PP^m\times\PP^s\rightarrow\PP^n\times\PP^m\times\PP^s$
formed by multihomogeneous polynomials of some multidegree.
\begin{prop}
The number of fixed points of a generic map $F:\PP^r\rightarrow\PP^r$ of degree $d$ is
$$1+d+\ldots+d^r.$$
\end{prop}
\begin{dem}
The following proof is standard in intersection theory.
The fixed points of a map $F:\PP^r\rightarrow\PP^r$ may be computed in $A(\PP^r\times\PP^r)$
as the degree of the product of the class of the graph of $F$, $[\Gamma]$, and the class of the diagonal, $[\triangle]$.
First, let's find out the class of the diagonal,
$$[\triangle]\in A^r(\PP^r\times\PP^r)=\Z[a,\alpha]/(a^{r+1},\alpha^{r+1}).$$
Being of codimension $r$, the class is a homogeneous polynomial of degree $r$,
$$[\triangle]=t_0\alpha^{r}+t_1a\alpha^{r-1}+\ldots+t_{r-1}a^{r-1}\alpha+t_ra^r,\quad t_i\in\Z$$
Here, $a$ represents a class of a hyperplane in $\PP^r$ and $a^i$
represents the intersection of $i$ of these generic hyperplanes, in other words, $a^i$ is a generic space
$\PP^{r-i}$ inside $\PP^r$. Same for $\alpha$ and $\alpha^j$.
Viewed in $\PP^r\times \PP^r$, $a^i$ is the class of $U\times\PP^r$,
$a^i=[U\times\PP^r]$,
and $\alpha^j$ is the class of $\PP^r\times V$, $\alpha^j=[\PP^r\times V]$, where $\dim U=r-i$ and $\dim V=r-j$.
The class $a^i\alpha^j$, represents a product of general linear spaces $U\times V\subseteq\PP^r\times\PP^r$,
where $\dim U=r-i$ and $\dim V=r-j$.\\

The class of the diagonal is determined by the coefficients $t_0,\ldots,t_r$.
Note that $t_i=[\triangle].a^{r-i}\alpha^{i}$.
Then, we need to count the number of points in $(U\times V)\cap \triangle$,
$$(U\times V)\cap \triangle\cong U\cap V=\{p\}\then $$
$$t_0=\ldots=t_r=1\then[\triangle]=\sum_{i=0}^ra^i\alpha^{r-i}.$$
Now, let's compute the class of the graph of a map, $\Gamma=\{(x,F(x))\}\subseteq\PP^r\times\PP^r$,
$$[\Gamma]\in A^r(\PP^r\times\PP^r)=\Z[a,\alpha]/(a^{r+1},\alpha^{r+1}),$$
it is also a homogeneous polynomial of degree $r$,
$$[\Gamma]=\tau_0\alpha^{r} +\tau_1a\alpha^{r-1}+\ldots+\tau_{r-1}a^{r-1}\alpha+\tau_ra^r,\quad \tau_i\in\Z.$$
Again, we have $\tau_i=[\Gamma].a^{r-i}\alpha^{i}$ so we need to count the points in $\Gamma\cap (U\times V)$,
where $\dim U=i$ and $\dim V=r-i$,
$$\Gamma\cap (U\times V)\cong\{x\in U\,|\,F(x)\in V\}=U\cap F^{-1}(V)\subseteq\PP^r.$$
If $F$ is formed by homogeneous polynomials of degree $d$, the pull-back of a hyperplane is a hypersurface of degree $d$, then
$$[U\cap F^{-1}(V)]=\alpha^i.F^\star(\alpha^{r-i})=\alpha^i.F^\star(\alpha)^{r-i}=\alpha^i(d\alpha)^{r-i}=d^{r-i}\alpha^r.$$
Then, $U\cap F^{-1}(V)$ has $d^{r-i}$ points, i.e. $\tau_i=d^{r-i}$,
$$[\Gamma]=d^r\alpha^{r} +d^{r-1}a\alpha^{r-1}+\ldots+da^{r-1}\alpha+a^r\then$$
$$[\triangle].[\Gamma]=(\sum_{i=0}^ra^i\alpha^{r-i})(\sum_{j=0}^rd^{r-j}a^j\alpha^{r-j})=
\sum_{i,j=0}^rd^{r-j}a^{i+j}\alpha^{2r-(i+j)}=\sum_{j=0}^rd^{r-j}=1+d+\ldots+d^r.$$
Given that a constant map has one fixed point, we use the convention $d^0=1$ for $d=0$.\qed
\end{dem}

Let's adapt the previous calculation to $\PP^n\times\PP^m\times\PP^s$.
We will compute the class of the diagonal and the class of a
graph, and then we will multiply them to obtain the number of fixed points.

\begin{teo}
The number of fixed points of a map $F=(F_1,F_2,F_3):\PP^n\times\PP^m\times\PP^s\rightarrow\PP^n\times\PP^m\times\PP^s$
is the coefficient of $\alpha^n\beta^m\gamma^s$ in the following polynomial in $\Z[\alpha,\beta,\gamma]$,
$$\sum_{i=0}^n\sum_{j=0}^m\sum_{k=0}^s(d_1\alpha+d_2\beta+d_3\gamma)^{n-i}(e_1\alpha+e_2\beta+e_3\gamma)^{m-j}
(f_1\alpha+f_2\beta+f_3\gamma)^{s-k}\alpha^{i}\beta^{j}\gamma^{k},$$
where $(d_1,d_2,d_3)$, $(e_1,e_2,e_3)$ and $(f_1,f_2,f_3)$ is the multidegree of $F_1$, $F_2$ and $F_3$ respectively.

For a generic map $F:\PP^{n_1}\times \ldots\times \PP^{n_k}\rightarrow \PP^{n_1}\times \ldots\times \PP^{n_k}$
the result is similar.
\end{teo}
\begin{dem}
The class of the diagonal $\triangle=\{(x,y,z,x,y,z)\}\in\PP^n\times\PP^m\times\PP^s\times \PP^n\times\PP^m\times\PP^s$,
is a homogenous polynomial of degree $n+m+s$,
$$[\triangle]\in
A(\PP^n\times\PP^m\times\PP^s\times \PP^n\times\PP^m\times\PP^s)=
\Z[\alpha,\beta,\gamma,a,b,c]/(\alpha^{n+1},\beta^{m+1},\gamma^{s+1},a^{n+1},b^{m+1},c^{s+1}),$$
Instead of doing the same computation as before, let
$$\pi_{1,4}:\PP^n\times\PP^m\times\PP^s\times \PP^n\times\PP^m\times\PP^s\rightarrow\PP^n\times\PP^n,$$
be the projection in the first and fourth factor (same for $\pi_{2,5}$ and $\pi_{3,6}$)
and let $\triangle_n\subseteq\PP^n\times\PP^n$ be the diagonal of $\PP^n$ (same for $\triangle_m$ and $\triangle_s$).
Then we have
$$[\triangle]=\pi_{1,3}^\star([\triangle_n]).\pi_{2,5}^\star([\triangle_m]).\pi_{3,6}^\star([\triangle_s])=
\sum_{i=0}^{n}\sum_{j=0}^{m}\sum_{k=0}^{s}a^i\alpha^{n-i}b^j\beta^{m-j}c^k\gamma^{s-k}.$$

The class of $\Gamma=\{(x,y,z,F(x,y,z))\}\subseteq\PP^n\times\PP^m\times\PP^s\times \PP^n\times\PP^m\times\PP^s$,
is a homogeneous polynomial of degree $n+m+s$,
$$[\Gamma]=\sum_{i+j+k+i'+j'+k'=n+m+s}\tau_{ijki'j'k'}a^i\alpha^{i'}b^j\beta^{j'}c^k\gamma^{k'},\quad \tau_{ijki'j'k'}\in\Z\then$$
$$\text{deg}([\triangle].[\Gamma])=\sum_{i=0}^{n}\sum_{j=0}^{m}\sum_{k=0}^{s}\tau_{ijkijk}.$$
Where \emph{deg} is the coefficient of $a^{n}\alpha^nb^m\beta^mc^s\gamma^s$.
Note that the integer $\tau_{ijkijk}$ is the degree
of $[\Gamma].a^{n-i}\alpha^ib^{m-j}\beta^jc^{s-k}\gamma^k$; the number of points in
$$\Gamma\cap (U_1\times U_2\times U_3\times V_1\times V_2\times V_3),$$
$$U_1,V_1\subseteq\PP^n,\,U_2,V_2\subseteq\PP^m,\,U_3,V_3\subseteq\PP^s,$$
$$\dim U_1+\dim V_1=n,\,\dim U_2+\dim V_2=m,\,\dim U_3+\dim V_3=s\then$$
$$\Gamma\cap (U_1\times U_2\times U_3\times V_1\times V_2\times V_3)\cong
(U_1\times U_2\times U_3)\cap F^{-1}(V_1\times V_2\times V_3)\subseteq\PP^n\times\PP^m\times\PP^s.$$
Let's use the fact that $F$ is equal to $(F_1,F_2,F_3)$,
$$F_1:\PP^n\times\PP^m\times\PP^s\rightarrow\PP^n,\,
F_2:\PP^n\times\PP^m\times\PP^s\rightarrow\PP^m,\,
F_3:\PP^n\times\PP^m\times\PP^s\rightarrow\PP^s$$
where $(d_1,d_2,d_3)$, $(e_1,e_2,e_3)$ and $(f_1,f_2,f_3)$ is the
multidegree of $F_1,F_2$ and $F_3$ respectively. Then
$$F^{-1}(V_1\times V_2\times V_3)=F_1^{-1}(V_1)\cap F_2^{-1}(V_2)\cap F_3^{-1}(V_3).$$
Thus, the class of the intersection that defines $\tau_{ijkijk}$ in the Chow ring $A(\PP^n\times\PP^m\times\PP^s)$, is
$$\tau_{ijkijk}=\alpha^{i}\beta^{j}\gamma^{k}F^\star(\alpha^{n-i}\beta^{m-j}\gamma^{s-k})=
\alpha^{i}\beta^{j}\gamma^{k}F_1^\star(\alpha^{n-i})F_2^\star(\beta^{m-j})F_3^\star(\gamma^{s-k})=$$
$$\alpha^{i}\beta^{j}\gamma^{k}(d_1\alpha+d_2\beta+d_3\gamma)^{n-i}(e_1\alpha+e_2\beta+e_3\gamma)^{m-j}(f_1\alpha+f_2\beta+f_3\gamma)^{s-k}.$$
\qed
\end{dem}

\begin{ejemplo}
Let's apply the previous formula to $\nabla\ell$
where $\ell:\mathbb{S}^{2}\times\mathbb{S}^{2}\times\mathbb{S}^{2}\rightarrow\R$ is a generic trilinear form.
The multidegree of $\partial\ell/\partial x$,
$\partial\ell/\partial y$ and $\partial\ell/\partial z$ is $(0,1,1)$, $(1,0,1)$ and $(1,1,0)$ respectively.
Then the number of fixed points of this map (over $\C$) is equal to $37$.
The number $37$, according to \cite[3]{1574201}, is the number of generalized singular values of $\ell$.
\end{ejemplo}

\begin{ejemplo}
Let's apply the formula to count the number of eigenvectors of a generic linear map $L:\R^{n+1}\rightarrow\R^{n+1}$.
The map $L$ induces a map $\PP^n\rightarrow\PP^n$ of degree $1$, then
$$\sum_{i=0}^n\alpha^{n-i}\alpha^i=\sum_{i=0}^n\alpha^{n}=(n+1)\alpha^n.$$
The map $\PP^n\rightarrow\PP^n$ has $n+1$ fixed points over $\C$, that is, $L$ has $n+1$ eigenvectors over $\C$.
\end{ejemplo}

\begin{ejemplo}
Finally, let's apply the formula to find the number of singular values of a generic linear map $L:\R^{n+1}\rightarrow\R^{m+1}$.
The map $L$ induces a bilinear form $\ell:\R^{n+1}\times \R^{m+1}\rightarrow\R$
and the bidegree of $\partial\ell/\partial x$ and $\partial\ell/\partial y$ is $(0,1)$ and $(1,0)$ respectively (assume $n\geq m$).
$$\sum_{i=0}^n\sum_{j=0}^m \beta^{n-i}\alpha^{m-j}\alpha^i\beta^{j}=
\sum_{i=0}^n\sum_{j=0}^m \beta^{n-i+j}\alpha^{m-j+i}=\sum_{j=0}^m \beta^{n-(n-m+j)+j}\alpha^{m-j+(n-m+j)}=(m+1)\alpha^n\beta^m.$$
Then, $\nabla\ell:\PP^n\times\PP^m\rightarrow\PP^n\times\PP^m$
has $m+1$ fixed points, that is, $L$ has $m+1$ singular values over $\C$.
We used the variational definition of singular values, see \cite{1574201}.
In case $n<m$ we can use the fact that the number of non-zero singular values of $L$ and $L^t:\R^{m+1}\rightarrow\R^{n+1}$
are the same.
\end{ejemplo}

\section{Theory for a bilinear form.}\label{sec-bil}
In this section we present a method to find the maximum of a bilinear form, $\ell$, over a product
of spheres, $\mathbb{S}^n\times \mathbb{S}^m$. This case is very special and the method presented here
does not work for a general multilinear form.

The key point of this method is the fact that the partial derivatives of $\ell=\sum a_{ij}x_iy_j$ are linear,
$$\frac{\partial\ell}{\partial x_i}(x,y)=\ell(e_i,y),\,\frac{\partial\ell}{\partial y_j}(x,y)=\ell(x,e_j),
\quad (x,y)\in\R^{n+1}\times\R^{m+1},\,0\leq i\leq n,\,0\leq j\leq m,$$
where $e_i$ and $e_j$ are canonical basis vectors of $\R^{n+1}$ and $\R^{m+1}$ respectively.
The map $\nabla\ell$ induces a linear map $L:\PP^{n+m+1}\rightarrow\PP^{n+m+1}$.
Let $(x_0:\ldots:x_n:y_0:\ldots:y_m)$ be a point in $\PP^{n+m+1}$. Then
$$L(x_0:\ldots:x_n:y_0:\ldots:y_m)=
\left(\ell(e_0,y):\ldots:\ell(e_n,y):\ell(x,e_0):\ldots:\ell(x,e_m)\right).$$
This map is well-defined. Let $\lambda\in\R$, $\lambda\neq 0$,
$$L(\lambda x_0:\ldots:\lambda x_n:\lambda y_0:\ldots:\lambda y_m)=
\left(\ell(e_0,\lambda y):\ldots:\ell(e_n,\lambda y):\ell(\lambda x,e_0):\ldots:\ell(\lambda x,e_m)\right)=$$
$$\left(\lambda\ell(e_0,y):\ldots:\lambda\ell(e_n,y):\lambda\ell(x,e_0):\ldots:\lambda\ell(x,e_m)\right)=
L(x_0:\ldots:x_n:y_0:\ldots:y_m).$$

\begin{teo}\label{bil-pn}
Let $p=(x,y)\in\mathbb{S}^n\times\mathbb{S}^m$ be an absolute maximum of $\ell$.
Then
$$\lim_{r\rightarrow\infty} L^r(q)=[p]$$
for a generic $q\in\PP^{n+m+1}$.
\end{teo}
\begin{dem}
Let $A\in\R^{n+m+2\times n+m+2}$ be a matrix representing the linear map $L$.
Given that $L$ is linear, the differential of $L$ at any point, $q$, is equal to $L$,
$$dL_q=L,\quad\forall q\in\PP^{n+m+1}.$$
In particular, the matrix $A$, also represents the differential of $L$ at $p$,
$$A=\widehat{dL_p}.$$
Let $\{v_0,\ldots,v_{n+m+1}\}$ be a basis of $\R^{n+m+2}$ formed by eigenvectors of $A$.
Let $\lambda_i$ be the eigenvalue of $v_i$, $0\leq i\leq n+m+1$.
By \autoref{def-euler} we know that $p$ is an eigenvector of $A$ with eigenvalue $\ell(p)$.
In particular, if the magnitude of $\lambda_0$ is maximum, then $|\lambda_0|=|\ell(p)|$ and $[p]=[v_0]$.

Let $z\in\R^{n+m+2}$ be a vector representing $q$ such that $z=a_0v_0+\ldots+a_{n+m+1} v_{n+m+1}$, $a_0\neq 0$.
Then
$$L^r(q)=[A^r.\sum_{i=0}^{n+m+1} a_iv_i]=[a_0\lambda_0^rv_0+\sum_{i=1}^{n+m+1} a_i\lambda_i^r v_i]=
[a_0v_0+\sum_{i=1}^{n+m+1} a_i\frac{\lambda_i^r}{\lambda_0^r} v_i]\rightarrow [a_0v_0]=p.$$
\qed
\end{dem}

In the proof of the previous theorem, we saw that the iterations of a linear map
in projective space converges to an eigenvector of maximum eigenvalue.
In particular, given a square matrix $A\in \R^{n+1\times n+1}$ and a generic vector $w\in\R^{n+1}$,
the sequence $\{[w],[Aw],[A^2w]\ldots\}\subseteq\PP^n$,
converges to a point $[v]$. The vector $v$ satisfies $Av=\lambda v$, where $|\lambda|$ is spectral radius of $A$.

The rate of convergence of this method is linear.

\begin{obs}\label{bil-ps1}
Let's give an algorithm to find the absolute maximum of a generic bilinear form,
$$\ell:\R^{n+1}\times\R^{m+1}\fl{}\R,$$
based on \autoref{bil-pn}.
Let $\nabla\ell=(\partial\ell/\partial x,\partial\ell/\partial y)$
be the gradient of $\ell$ and let $q$ be the initial condition, where $q=(x,y)$, $x\in\R^{n+1}$, $x\neq0$
and $y\in\R^{m+1}$, $y\neq0$.

{\footnotesize
$$\begin{array}{rl}
\text{\tt Input:}&\text{\tt A bilinear form }\ell:\R^{n+1}\times\R^{m+1}\rightarrow\R.\\
\text{\tt Output:}&\text{\tt The absolute maximum }(x,y)\in\mathbb{S}^{n}\times\mathbb{S}^{m}.\\
\hline
1.&\text{\tt Let }q=q/\|q\|\,\,\text{\tt and }aux=(1,0,\ldots,0).\\
2.&\text{\tt While }|\langle q,aux\rangle|\,\,\text{\tt is different from }1\,\,\text{\tt do}\\
&\begin{array}{rl}
2.1& aux=q\\
2.2& q=\nabla\ell(q)\\
2.3& q=q/\|q\|\\
\end{array}\\
3.&\text{\tt Let }x=(q_0,\ldots,q_n),\,y=(q_{n+1},\ldots,q_{n+m+1}).\\
4.&\text{\tt Return } (x/\|x\|,y/\|y\|).
\end{array}$$
}The iterations stops when the points in projective space are equal, in other words, when the
cosine of the angle between $q$ and $aux$ is $1$ or $-1$ (when they are aligned).
Given that the absolute maximum is attractive (see \autoref{bil-pn}), the program ends.
The maximum value is $|\ell(x,y)|$.
\end{obs}

\begin{obs}
We may adapt this algorithm to a multilinear form, but in the multilinear case, in general,
the absolute maximum is not an attractive fixed point.
For example, the trilinear form $\ell:\R^2\times\R^2\times\R^2\rightarrow\R$,
{\small
$$\ell(x,y,z)=6x_0y_0z_0+3x_1y_0z_0-6x_0y_1z_0+16x_1y_1z_0-14x_0y_0z_1-15x_1y_0z_1-11x_0y_1z_1+8x_1y_1z_1,$$
}induces a map $\PP^{5}\rightarrow\PP^5$ of degree $2$ without attractive fixed points.
Even more, the $4$-multilinear form $\ell:\R^2\times\R^2\times\R^2\times\R^2\rightarrow\R$,
{\small
$$\begin{array}{rcl}
\ell(x,y,z,v)&=&
4x_0y_0z_0v_0+6x_1y_0z_0v_0+x_0y_1z_0v_0-6x_1y_0z_1v_0v_0-5x_0y_0z_1v_0+\\
& &7x_1y_1z_0-5x_0y_1z_1v_0+2x_0y_0z_0v_1-3x_1y_0z_0v_1-7x_0y_1z_0v_1+\\
& &9x_1y_1z_0v_1-9x_0y_0z_1v_1-9x_1y_0z_1v_1-6x_0y_1z_1v_1+8x_1y_1z_1v_1,
\end{array}$$
}induces a map $\PP^7\rightarrow\PP^7$ of degree $3$ with
two attractive fixed points. One is the absolute maximum.
\end{obs}

\section{Presentation of the general algorithm.}\label{sec-pres}
In this section we present an algorithm to find the maximum of a multilinear form over a product of
spheres. First, we reduce the problem to a system of multilinear equations and then
we resolve the system using algebro-geometric tools.
In the first part of the section, we present an algorithm to find the absolute maximum of a
multilinear form. In the second, we give an algorithm to find the point where the maximum occurs.
This last algorithm requires some extra hypothesis.

\begin{prop}\label{pres-sys}
Let $\ell:\R^{n+1}\times\R^{m+1}\times\R^{s+1}\rightarrow\R$ be a generic trilinear form.
There exists a bijection between classes of extreme points of $\ell$
and solutions of the following system of trilinear equations in $\PP^n\times\PP^m\times\PP^s$,
$$\left\{\begin{array}{rcl}
\ell(x_je_i-x_ie_j,y,z)&=&0,\quad 0\leq i<j\leq n,\\
\ell(x,y_je_i-y_ie_j,z)&=&0,\quad 0\leq i<j\leq m,\\
\ell(x,y,z_je_i-z_ie_j)&=&0,\quad 0\leq i<j\leq s,\\
\end{array}\right.$$
The vector $e_k$ satisfies $(e_k)_l=0$ if $l\neq k$ and $(e_k)_k=1$.

In the multilinear case, we obtain a similar result; a system of multilinear equations.
\end{prop}
\begin{dem}
From \autoref{theo-bij} we know that every class of an extreme point of $\ell$, is a fixed point of
$\nabla\ell:\PP^n\times\PP^m\times\PP^s\rightarrow\PP^n\times\PP^m\times\PP^s$.
If $\ell$ is a generic trilinear form, we know that the number of fixed points is finite (see \autoref{sec-num}).

A fixed point of $\nabla\ell$, $([x],[y],[z])$, satisfies
$$\left\{\begin{array}{lcr}
\partial\ell/\partial x_i(x,y,z)&=&2\alpha x_i,\quad 0\leq i\leq n,\\
\partial\ell/\partial y_j(x,y,z)&=&2\beta y_j,\quad 0\leq j\leq m,\\
\partial\ell/\partial z_k(x,y,z)&=&2\lambda z_k,\quad 0\leq k\leq s,
\end{array}\right.$$
where $\alpha$, $\beta$ and $\lambda$ are three nonzero real numbers.
In $\PP^n\times\PP^m\times\PP^s$, the equations are
$$\left\{\begin{array}{lcr}
x_j\partial\ell/\partial x_i(x,y,z)&=& x_i\partial\ell/\partial x_j(x,y,z),\quad 0\leq i<j\leq n,\\
y_j\partial\ell/\partial y_i(x,y,z)&=& y_i\partial\ell/\partial y_j(x,y,z),\quad 0\leq i<j\leq m,\\
z_j\partial\ell/\partial z_i(x,y,z)&=& z_i\partial\ell/\partial z_j(x,y,z),\quad 0\leq i<j\leq s,
\end{array}\right.$$
The result follows from the equalities,
$$\partial\ell/\partial x_i(x,y,z)=\ell(e_i,y,z),\quad
\partial\ell/\partial y_j(x,y,z)=\ell(x,e_j,z),\quad
\partial\ell/\partial z_k(x,y,z)=\ell(x,y,e_k).$$
\qed
\end{dem}

Let's present the algorithm to find the absolute maximum of a generic multilinear form.
The resolution of the system is bases on \emph{Eigenvalue Theorem}. Let's recall it.
Consider a system of polynomial equations with finitely many solutions in $\C^{n}$,
$$\left\{\begin{array}{rcl}
f_1(x_1,\ldots,x_n)&=&0\\
&\vdots&\\
f_m(x_1,\ldots,x_n)&=&0\\
\end{array}\right.$$
where $f_1,\ldots,f_m$ are polynomials in $\C[x_1,\ldots,x_n]$. The quotient ring,
$$\mathcal{A}=\C[x_1,\ldots,x_n]/\langle f_1,\ldots,f_m\rangle,$$
is a finite-dimensional vector space, \cite[Theorem 2.1.2]{MR2161984}.
The dimension of $\mathcal{A}$ is the number of solutions of the system.

Every polynomial $f\in\C[x_1,\ldots,x_n]$, determines a linear map $M:\mathcal{A}\rightarrow \mathcal{A}$,
$$M(\overline{g})=\overline{fg},\quad g\in\C[x_1,\ldots,x_n],$$
where $\overline{g}$ denotes the class of the polynomial $g$ in the quotient ring $\mathcal{A}$.
The matrix of $M$ is called the \emph{multiplication matrix} associated to the polynomial $f$.
\begin{nteo}[Eigenvalue Theorem]
The eigenvalues of $M$ are $\{f(p_1),\ldots,f(p_r)\}$, where $\{p_1,\ldots,p_r\}$
are the solutions of the system of polynomial equations.
See \cite[Theorem 2.1.4]{MR2161984} for a proof.
\end{nteo}

The algorithm in \autoref{app-ps2}, first generates the following system of polynomial equations,
$$\left\{
\begin{array}{cc}
\begin{array}{lcr}
x_j\partial\ell/\partial x_i(x,y,z)&=& x_i\partial\ell/\partial x_j(x,y,z),\quad 0\leq i<j\leq n,\\
y_j\partial\ell/\partial y_i(x,y,z)&=& y_i\partial\ell/\partial y_j(x,y,z),\quad 0\leq i<j\leq m,\\
z_j\partial\ell/\partial z_i(x,y,z)&=& z_i\partial\ell/\partial z_j(x,y,z),\quad 0\leq i<j\leq s,\\
\end{array}\\
x_0^2+\ldots+x_n^2=1,\quad y_0^2+\ldots+y_m^2=1,\quad z_0^2+\ldots+z_s^2=1.
\end{array}
\right.$$
Then, computes the real eigenvalues, $\{\lambda_0,\ldots,\lambda_r\}$,
of the multiplication matrix associated to $\ell$.
Finally, it returns $\lambda_i$ such that $|\lambda_i|\geq |\lambda_j|$ for all $0\leq j\leq r$.
This number, is the maximum of $\ell$ over $\mathbb{S}^n\times\mathbb{S}^m\times\mathbb{S}^s$.

For the algorithm and an implementation in Maple, see \autoref{app-ps2}.\\

Now, let's give an algorithm
to find the point $(x,y,z)\in\mathbb{S}^n\times\mathbb{S}^m\times\mathbb{S}^s$ such that $|\ell(x,y,z)|$
is maximum. We need to use the following result (same notation as Eigenvalue Theorem),
\begin{nteo}
Let $x=\lambda_1x_1+\ldots+\lambda_nx_n$ be a generic linear form and let $M$ be its multiplication matrix.
Assume that $B=\{1,x_1,\ldots,x_n,\ldots\}$ is a finite basis of $\mathcal{A}$.
Then, the eigenvectors of $M$ determine solutions of the
system of polynomial equations. Specifically, if $v=(v_0,\ldots,v_n,\ldots)$
is an eigenvector of $M$ such that $v_0=1$, then $(v_1,\ldots,v_n)$ is a solution
of the system of polynomial equations. Even more, every solution is of this form.
See \cite[\S 2.1.3]{MR2161984} for a proof.
\end{nteo}
Note that the Theorem requires that the variables $\{x_1,\ldots,x_n\}$ are elements of the basis $B$.
It could be the case that some variables are missing from $B$.
For example, if
$x_1,\ldots,x_i\in B$, and $x_{i+1} ,\ldots,x_n\not\in B$,
then, every missing variable, say $x_j$, is a linear combination of $\{x_{1} ,\ldots,x_{i}\}$,
$$x_j=a_{j1}x_1+\ldots+a_{ji}x_i,\quad i+1\leq j\leq n$$
If $v=(1,v_1,v_2\ldots)$ is an eigenvector of $M$, the $j$-coordinate of the solution
corresponding to $v$, is $a_{j1}v_1+\ldots+a_{ji}v_i$. See \cite[\S 2.1.3]{MR2161984}.\\

In order to apply the previous Theorem, we need to guarantee that the basis $B$ contains all the variables.
The \emph{affine system} in \autoref{pres-sys} is,
$$\left\{
\begin{array}{cc}
\begin{array}{lcr}
x_j\partial\ell/\partial x_i(x,y,z)&=& x_i\partial\ell/\partial x_j(x,y,z),\quad 0\leq i<j\leq n,\\
y_j\partial\ell/\partial y_i(x,y,z)&=& y_i\partial\ell/\partial y_j(x,y,z),\quad 0\leq i<j\leq m,\\
z_j\partial\ell/\partial z_i(x,y,z)&=& z_i\partial\ell/\partial z_j(x,y,z),\quad 0\leq i<j\leq s,\\
\end{array}\\
x_0=1,\quad y_0=1,\quad z_0=1.
\end{array}
\right.$$
The solutions of this system determine classes of extreme points of $\ell$.
The genericity of $\ell$ implies that all the extreme points of $\ell$, $(x,y,z)$, satisfy $x_0\neq 0,y_0\neq 0,z_0\neq 0$.
Then, all the classes of extreme points appear as the solutions of the affine system in \autoref{pres-sys}.

\begin{teo}
Assume that $\ell:\R^{n+1}\times\R^{m+1}\times \R^{s+1}\rightarrow\R$ is
a generic trilinear form and that $2n,2m,2s\leq n+m+s$, see \autoref{theo-hyp}.
Then, the basis $B$ of the affine system in \autoref{pres-sys}
contains all the variables.

In the multilinear case, we obtain a similar result.
\end{teo}
\begin{dem}
Given that the equations in \autoref{pres-sys} are multilinear, the quotient ring, $\mathcal{A}$,
is multi-graded. 
Let's denote $\mathcal{A}_{(d_1,d_2,d_3)}$ the
multidegree part $(d_1,d_2,d_3)$, where $d_1,d_2,d_3\geq 0$.
The hypothesis $2n,2m,2s\leq n+m+s$, implies that the following set is empty,
$$\left\{(x,y,z)\in\PP^n\times\PP^m\times\PP^s\,|\,\frac{\partial\ell}{\partial x}(x,y,z)=0\text{ or }
\frac{\partial\ell}{\partial y}(x,y,z)=0\text{ or }\frac{\partial\ell}{\partial z}(x,y,z)=0\right\}=\emptyset.$$
Then, the equations $\{\partial\ell/\partial x_i\}_{i=0}^n$ are linearly independent. Same for 
$\{\partial\ell/\partial y_j\}_{j=0}^m$ and $\{\partial\ell/\partial z_k\}_{k=0}^s$.
In the quotient ring, $\mathcal{A}$, the partial derivatives, are proportional to
the variables, thus, the variables are linearly independent too.
For example, a basis for the multidegree part $(0,0,0)$ is $\{1\}$, and a basis
for the multidegree part $(1,0,0)$ is $\{x_0,\ldots,x_n\}$. Even more,
a basis for
$$\mathcal{A}_{(1,0,0)}\oplus\mathcal{A}_{(0,1,0)}\oplus\mathcal{A}_{(0,0,1)}$$
is $\{x_0,\ldots,x_n,y_0,\ldots,y_m,z_0,\ldots,z_s\}$.

Let's add the equations $x_0=y_0=z_0=1$ to the system of polynomial equations.
The equations are not multilinear, so the corresponding quotient ring is not multi-graded,
$$\widehat{\mathcal{A}}=\mathcal{A}/\langle x_0-1,y_0-1,z_0-1\rangle.$$

Let's see that the variables $\{x_1,\ldots,x_n,y_1,\ldots,y_m,z_1,\ldots,z_s\}$
are linearly independent in $\widehat{\mathcal{A}}$. This implies that the basis $B$
of $\widehat{\mathcal{A}}$, formed by monomials, contains all the variables.
$$\sum_{i=1}^n\alpha_ix_i+\sum_{j=1}^m\beta_jy_j+\sum_{k=1}^s\lambda_kz_k=0\in\widehat{\mathcal{A}},\quad
\alpha_i,\beta_j,\lambda_k\in\C\then$$
$$\sum_{i=1}^n\alpha_ix_i+\sum_{j=1}^m\beta_jy_j+\sum_{k=1}^s\lambda_kz_k=P\in\mathcal{A},\quad P\in\langle x_0-1,y_0-1,z_0-1\rangle,$$
where $P$ is a polynomial combination of $x_0-1$, $y_0-1$ and $z_0-1$.

Denote $P_{(d_1,d_2,d_3)}$ the multidegree part, $(d_1,d_2,d_3)$, of $P$.
Given that $\mathcal{A}$ is multi-graded, we get the following equalities in $\mathcal{A}$,
$$P=P_{(1,0,0)}+P_{(0,1,0)}+P_{(0,0,1)},\quad
\sum_{i=1}^n\alpha_ix_i=P_{(1,0,0)},\quad \sum_{j=1}^m\beta_jy_i=P_{(0,1,0)},\quad \sum_{k=1}^s\lambda_kz_k=P_{(0,0,1)}.$$
Using the fact that the variables $\{x_0,\ldots,x_n\}$ are linearly independent in $\mathcal{A}$,
we obtain that $x_0$ is not a variable in $P$. Same for $y_0$ and $z_0$.
Given that $P$ is a polynomial combination of $x_0-1$, $y_0-1$ and $z_0-1$, it must be $0$.
Then, $\alpha_1=\ldots=\alpha_n=0$, $\beta_1=\ldots=\beta_m=0$ and $\lambda_1=\ldots=\lambda_s=0$.
\qed
\end{dem}

\begin{obs}\label{pres-algo-b-2}
The algorithm to find the point $(x,y,z)\in\mathbb{S}^n\times\mathbb{S}^m\times\mathbb{S}^s$
such that $|\ell(x,y,z)|$ is maximum, is the following.

The reader may adapt the algorithm to a multilinear form.
{\footnotesize
$$\begin{array}{rl}
\text{\tt Input:}&\text{\tt A generic trilinear form }\ell:\R^{n+1}\times\R^{m+1}\times\R^{s+1}\rightarrow\R,\\
&\text{\tt where }2n+2m+2s\leq n+m+s.\\
\text{\tt Output:}&\text{\tt The absolute maximum }(x,y,z)\in\mathbb{S}^{n}\times\mathbb{S}^{m}\times\mathbb{S}^{s}.\\
\hline
1.&\text{\tt Compute the system of trilinear equations of \autoref{pres-sys}}.\\
2.&\text{\tt Add the equations }x_0=y_0=z_0=1.\\
3.&\text{\tt Compute a G\"obner basis for the resulting system, }\mathcal{I}.\\
4.&\text{\tt Find a basis }B\,\,\text{\tt of }\C[x_1,\ldots,x_{n},y_1,\ldots,y_m,z_1,\ldots,z_s]/\mathcal{I}.\\
5.&\text{\tt Compute the multiplication matrix of }x_1.\\
6.&\text{\tt Compute the eigenvectors of the multiplication matrix.}\\
7.&\text{\tt For each eigenvector  }v\,\,\text{\tt do }\\
&\begin{array}{rl}
7.1&\text{\tt Normalize }v\,\,\text{\tt such that }v=(1,v_1,\ldots).\\
7.2&\text{\tt Let }x=(x_1,\ldots,x_n)\,\,\text{\tt be such that }x_i=v_{\sigma_i}\,\,\text{\tt where}\\
&\sigma_i\,\,\text{\tt is the coordinate of }x_i\,\,\text{\tt in }B,\,1\leq i\leq n.\\
7.3&\text{\tt Same for }y\,\,\text{\tt and }z.\\
7.4&\text{\tt Normalize the points, }x=x/\|x\|,\,y=y/\|y\|,\,z=z/\|z\|.\\
7.5&\text{\tt Evaluate }\ell\,\,\text{\tt at }(x,y,z)\,\,\text{\tt if the coordinates are real.}\\
7.6&\text{\tt Save the maximum.}
\end{array}\\
8.&\text{\tt Return the maximum, }(x,y,z).\\
\end{array}$$
}

In Step 5 of the algorithm we used the linear form $x_1$ as a generic linear form.
This fact is not restrictive. Given that the trilinear form is generic,
we may suppose that the first coordinates of the classes of extreme points of $\ell$ are all different.
In other words, the eigenvalues of the multiplication matrix of $x_1$ have multiplicity one.
See \cite[\S 2.1.3]{MR2161984}.
\end{obs}

\section{Applications and examples.}\label{sec-ex}
Let's start we some applications. First, we give applications of the iterative algorithm
to find the maximum of a bilinear form. Then, we give applications of the general algorithm.

\begin{obs}
Given a real matrix $A$, its first singular value (the 2-norm) is given by
$$\max_{\|x\|=\|y\|=1} x^tAy.$$
In other words, it is the maximum over a product of spheres of the bilinear form $(x,y)\rightarrow x^t A y$.\\

An interesting fact of \autoref{bil-pn} is that we can find the first singular vectors and the first singular value, $|\ell(x,y)|$,
of $\ell$ without using the \emph{spectral radius formula}.
Recall that the 2-norm of a matrix $A$ is computed using the spectral radius formula,
$$\|A\|_2=\sqrt{\lim_{k\rightarrow+\infty}\|(AA^t)^k\|^{\frac{1}{k}}}.$$
\end{obs}

\begin{ejemplo}
Let $A\in\R^{4\times 3}$ be the matrix
$$A=\begin{pmatrix}
3&2&32\\
2&1&36\\
-3&25&2\\
0&-1&1\\
\end{pmatrix}$$
Then, with the algorithm in \autoref{bil-ps1},
we get that the 2-norm is $48.46054603$. Using the standard algorithm (the spectral radius of $AA^t$)
we get the same number $48.46054603$.
\end{ejemplo}

\begin{ejemplo}
An interesting example is the maximum of a bilinear form over $\mathbb{S}^1\times\mathbb{S}^0$.
Note that the domain is a cylinder of $\R^3$, so we can draw the whole situation.
Take, for example, the bilinear form
$$\ell:\R^2\times\R\fl{}\R,\quad\ell(x,y)=4x_1y+2x_2y.$$
The maximum of $\ell$ over $\mathbb{S}^1\times\mathbb{S}^0$ is the 2-norm of the vector $(4,2)$, that is,
$$\|(4,2)\|=\sqrt{20}\cong 4.472135954.$$
Let's compute this using the algorithm in \autoref{bil-ps1}. First of all note that the gradient of $\ell$
determines a vector field over the cylinder, and the iteration follows the arrow.
Over the ending point of the iteration, the flow is orthogonal to the surface.
This means, that we have reached an extreme,
$$\max_{\|x\|=|y|=1}\ell(x,y)\cong 4.472135953.$$
\end{ejemplo}

Let's give now some applications of the general algorithm, see \autoref{app-ps2}.
\begin{obs}\label{ex-quantum}
The first interesting application of the algorithm in \autoref{app-ps2} is to the theory of entanglement.
It is of interest to find
the maximum of the form $\langle\rho,-\rangle$ over the space of \emph{separable states}. The matrix
$\rho$ is called a \emph{state} if it is hermitian, $\rho\geq0$ and $\text{tr}(\rho)=1$.
It is easy to see that the space of states is a convex set and is generated by the
matrices of the form $\rho_i\rho_i^\dag$ where $\rho_i$ is a column vector of norm one in a finite dimensional
vector space $\mathcal{H}$, $\rho_i\in\mathcal{H}$, $\|\rho_i\|=1$.
The general theory says that when we are working with two particles, we need to consider the space
of states over the tensor product $\mathcal{H}=\mathcal{H}_1\otimes\mathcal{H}_2$. In this situation
a state is called \emph{separable} if
it is a convex combination of the form $\sum a_{i}v_i\otimes w_i$, where $v_i$ is a state of $\mathcal{H}_1$ and
$w_i$ a state of $\mathcal{H}_2$.
Let's call $Sep(\mathcal{H})$ the convex space of separable states.
It is true that the space of separable states is a convex set generated by the matrices of the form
$xx^\dag\otimes yy^\dag$, where $x\in\mathcal{H}_1$, $y\in\mathcal{H}_2$ and $\|x\|=\|y\|=1$. Then
$$\max_{Sep(\mathcal{H})}\langle\rho,-\rangle=\max_{\|x\|=\|y\|=1}\langle\rho,xx^\dag\otimes yy^\dag\rangle.$$
Note that the form is not bilinear in $x$ nor in $y$.
Rewriting the state $\rho$ in the form $\rho=\sum\lambda_i\rho_i\rho_i^\dag$ with $\langle\rho_i,\rho_j\rangle=0$, $\|\rho_i\|=1$,
and using the equality $\langle \rho_i,x\otimes y\rangle^2=\langle \rho_i\rho_i^\dag,xx^\dag\otimes yy^\dag\rangle$, we get
$$\max_{Sep(\mathcal{H})}\langle\rho,-\rangle=
\max_{\|x\|=\|y\|=1}\|\sum \sqrt{\lambda_i}\langle \rho_i,x\otimes y\rangle \rho_i\|^2.$$
The resulting map $\sum \sqrt{\lambda_i}\langle \rho_i,x\otimes y\rangle \rho_i$
is bilinear in $x$ and in $y$ and our algorithm is capable to maximize it.
See the next example.
\end{obs}

\begin{ejemplo}
Suppose that we are working with the following state $\rho$ in $\R^4=\R^2\otimes\R^2$,
{\scriptsize
$$\begin{pmatrix}
0.242894940524649938&-0.123994312358229969&-0.0712215842649899789& 0.219784373378769966\\
-0.123994312358229969& 0.0888784895376599772&0.111143109132249979&-0.0627261109839499926\\
-0.0712215842649899789& 0.111143109132249979& 0.361255602168969903&0.0603142605185699871\\
0.219784373378769966&- 0.0627261109839499926& 0.0603142605185699871& 0.306970967813849916
\end{pmatrix}.$$
}
We choose to work over the real numbers to make the exposition clear, but all the results
can be adapted to work with hermitian matrices instead of symmetric matrices.
Using Cholesky and the Singular Value Decomposition Algorithm we have $\rho=\sum\lambda_i\rho_i\rho_i^\dag$,
{\footnotesize
$$\lambda_1=0.5435016101,\,\lambda_2=0.4146107959,$$
$$\lambda_3=0.04113792919,\,\lambda_4=0.0007496649711,\quad\lambda_1+\lambda_2+\lambda_3+\lambda_4=1.$$
$$\rho_1=\begin{pmatrix}
-0.656481390369177854\\
0.326643787198963642\\
0.245965753592146592\\
-0.633906040705653040
\end{pmatrix},\,
\rho_2=\begin{pmatrix}
-0.0253829550629408562\\
-0.209402292907094082\\
-0.881013881627254691\\
-0.423463015737623016
\end{pmatrix}$$
$$\rho_3=\begin{pmatrix}
-0.444223726945872255\\
0.546710519336902400\\
-0.399219012690601172\\
0.586853532298337144
\end{pmatrix},\,
\rho_4=\begin{pmatrix}
-0.609141338368644146\\
-0.741998735885678107\\
0.0627659808038323054\\
0.272846362424434330
\end{pmatrix}.$$
}
With the algorithm in \autoref{app-ps2} to the trilinear form
$$\ell(x,y,z)=\sum\sqrt{\lambda_i}\langle \rho_i,x\otimes y\rangle\langle \rho_i,z\rangle,\quad
(x,y,z)\in\mathbb{S}^1\times\mathbb{S}^1\times\mathbb{S}^3,$$
we get that the maximum is $0.7228016991$.
Finally
$$\max_{Sep(\mathcal{H})}\langle\rho,-\rangle=
\max_{\|x\|=\|y\|=1}\|\sum \sqrt{\lambda_i}\langle \rho_i,x\otimes y\rangle \rho_i\|^2=
\max_{\|x\|=\|y\|=\|z\|=1}|\ell(x,y,z)|^2\cong$$
$$0.7228016991^2\cong0.5224422962.$$
Note that if $\rho$ is separable then $\langle\rho,\rho\rangle\leq \max_{Sep(\mathcal{H})}\langle\rho,-\rangle$. This is
not the case, but for example, the following state is not separable (it is called \emph{entangled state}),
{\scriptsize
$$\begin{pmatrix}
0.168106937369559950&-0.190509527669719958&-0.200004375511779936&-0.0690454833860399825\\
-0.190509527669719958&0.257651665981429912&0.267759084652009926&0.0985801483325399742\\
-0.200004375511779936&0.267759084652009926&0.320790216378169901&0.194053687463299957\\
-0.0690454833860399825&0.0985801483325399742&0.194053687463299957&0.253451180300149959
\end{pmatrix}.$$
}
We have $\langle\rho,\rho\rangle\cong 0.6620536187\not\leq 0.4862909489\cong\max_{Sep(\mathcal{H})}\langle\rho,-\rangle$.
\end{ejemplo}

\begin{obs}\label{ex-rk-1}
Our final application is the ability to find numerically the closest rank-one tensor of a given tensor.
In the article \cite{MR2447444}, the authors considered the problem of finding the best rank-$r$ approximation
of a given tensor. They proved that for $r>1$ the problem is ill-posed, but when $r=1$ the problem has a solution, \cite[4.5]{MR2447444}.
Here we find the solution.
Let's prove that a computation of the absolute maximum of $\ell$ over a product of spheres
gives the closest \emph{rank-one multilinear form} to $\ell$.
A rank-one multilinear form is a product of linear forms, $\ell_1\ldots\ell_s$,
where $\ell_i:\R^{n_i+1}\rightarrow\R$, $1\leq i\leq s$.
We choose to do this remark about multilinear forms, but dually, the same is true for tensors.\\

For simplicity, we do the proof for a trilinear form. The proof is similar in the multilinear case.
Consider the affine Segre map (it is not an isometry)
$$\R^{n+1}\times \R^{m+1}\times \R^{s+1}\fl{}\R^{n+1}\otimes \R^{m+1}\otimes \R^{s+1},\quad (x,y,z)\fl{}x\otimes y\otimes z.$$
Using the usual inner product in the tensor product, we identify
$$\R^{n+1}\otimes \R^{m+1}\otimes \R^{s+1}\cong
(\R^{n+1}\otimes \R^{m+1}\otimes \R^{s+1})^\vee,$$
$$x\otimes y\otimes z\fl{}\langle x\otimes y\otimes z,-\rangle,\quad
\langle x\otimes y\otimes z,a\otimes b\otimes c\rangle=
\langle x,a\rangle\langle y,b\rangle\langle z,c\rangle.$$
We can identify the following three different notations
$$\ell(x,y,z)=\ell(x\otimes y\otimes z)=\langle \ell,x\otimes y\otimes z\rangle.$$
The first equality identifies a trilinear form with a linear map $\ell:\R^{n+1}\otimes \R^{m+1}\otimes \R^{s+1}\rightarrow\R$.
The second equality identifies, under the isometry
$(\R^{n+1}\otimes \R^{m+1}\otimes \R^{s+1})^\vee\cong\R^{n+1}\otimes \R^{m+1}\otimes \R^{s+1}$,
the linear form $\ell$ with the tensor $\ell\in\R^{n+1}\otimes \R^{m+1}\otimes \R^{s+1}$.

Let $\mathbb{S}$ be the immersion of $\mathbb{S}^n\times\mathbb{S}^m\times\mathbb{S}^s$ under the Segre map,
$$\mathbb{S}=\left\{\langle x\otimes y\otimes z,-\rangle\colon\,\|x\|=\|y\|=\|z\|=1\right\}
\subseteq(\R^{n+1}\otimes \R^{m+1}\otimes \R^{s+1})^\vee.$$
Then, for all
$\phi=\langle x\otimes y\otimes z,-\rangle\in \mathbb{S}$, we have
$$\|\ell-\phi\|^2=\left(\|\ell\|^2+\|\phi\|^2-2\langle\ell,\phi\rangle\right)=
\|\ell\|^2+1-2\ell(x,y,z).$$
In other words, a local maximum of $\ell$ is a local minimum of the distance function, $\|\ell-\phi\|$.

Let $\mathbb{B}$ be the image, under the Segre map, of a product of balls,
$$\mathbb{B}=\{\langle x\otimes y\otimes z,-\rangle\colon\,\|x\|,\|y\|,\|z\|\leq1\}\subseteq(\R^{n+1}\otimes \R^{m+1}\otimes \R^{s+1})^\vee.$$
Note that the elements of $\mathbb{B}$ are rank-one multilinear forms.
It is easy to see that $\mathbb{B}$ is compact and convex,
so the distance from $\ell$ to $\mathbb{B}$ is achieved in $\mathbb{S}$ (the border).
In other words, the closest rank-one multilinear form to $\ell$ is an element of $\mathbb{S}$.
Summing up, a computation with the algorithm in \autoref{pres-algo-b-2} of the absolute maximum of $\ell$,
gives a closest rank one multilinear form to $\ell$.
\end{obs}

\begin{ejemplo}
Let $\ell:\R^2\times\R^2\times\R^2\times\R^2\fl{}\R$ be a multilinear form,
$$\ell(x,y,z,t)=4x_1y_1z_1t_1+6x_2y_1z_1t_1+x_1y_2z_1t_1+7x_2y_2z_1t_1-5x_1y_1z_2t_1-$$
$$6x_2y_1z_2t_1-5x_1y_2z_2t_1+2x_1y_1z_1t_2-3x_2y_1z_1t_2-7x_1y_2z_1t_2+9x_2y_2z_1t_2-$$
$$9x_1y_1z_2t_2-9x_2y_1z_2t_2-6x_1y_2z_2t_2+8x_2y_2z_2t_2.$$
Using the algorithm in \autoref{pres-algo-b-2}, we get that the closest rank one multilinear form is
$$\ell_1(x)\ell_2(y)\ell_3(z)\ell_4(t),$$
$$\begin{array}{lcr}
\ell_1(x)&=&0.4799354720x_1-0.8773037918x_2\\
\ell_2(y)&=&0.2732019392y_1-0.9619567040y_2\\
\ell_3(z)&=&0.7563638894z_1+0.6541511043z_2\\
\ell_4(t)&=&0.3260948315t_1+0.9453370622t_2.
\end{array}$$
The value of the absolute maximum of $\ell$ is $16.71262553$.
\end{ejemplo}

\begin{ejemplo}
Let $v\in\R^2\otimes\R^3$ be the following tensor
$$v=4x_1\otimes y_1-9x_2\otimes y_1+2x_1\otimes y_2+x_2\otimes y_2-5x_1\otimes y_3-7x_2\otimes y_3.$$
Using the algorithm in \autoref{pres-algo-b-2} we get that the closest rank one tensor is
{\small
$$(0.01162554952x_1+0.9999324213x_2)\otimes(-0.7821828869y_1+0.08939199251y_2-0.6166027924y_3)$$
}In this case, we can check this result. The first singular vectors of the matrix
$$\begin{pmatrix}
4&-9\\
2&1\\
-5&-7
\end{pmatrix},$$
are
$$(0.01162554952, 0.99993242102),\quad(-0.7821828866, 0.08939199251,- 0.6166027924).$$
\end{ejemplo}

\appendix

\section{General Algorithm for a multilinear form. }\label{app-ps2}
Let's give an algorithm to find the maximum value of generic multilinear map over a product of spheres,
$$\ell:\R^{n_1}\times\ldots\times\R^{n_r}\rightarrow\R^{n_{r+1}},\quad \max_{\|x_1\|=\ldots=\|x_r\|=1}\|\ell(x_1,\ldots,x_r)\|.$$
Recall from \autoref{sec-theo} that we may assume that $\ell$ is a multilinear form,
$$\widehat{\ell}:\R^{n_1}\times\ldots\times\R^{n_r}\times\R^{n_{r+1}}\rightarrow\R
,\quad \max_{\|x_1\|=\ldots=\|x_{r+1}\|=1}|\widehat{\ell}(x_1,\ldots,x_{r+1})|.$$

The following is a pseudocode in the trilinear case.
We choose to work with Gr\"obner Basis because it is implemented in
most computer algebra systems (Maple, Macaulay2, Singular).
In \cite{MR1022943}, the authors proposed an algorithm without the need of Gr\"obner Basis.
See also \cite[2.3.1]{MR2161984}.

{\footnotesize
$$\begin{array}{rl}
\text{\tt Input:}&\text{\tt A generic trilinear form }\ell:\R^{n+1}\times\R^{m+1}\times\R^{s+1}\rightarrow\R.\\
\text{\tt Output:}&\text{\tt The absolute maximum }(x,y,z)\in\mathbb{S}^{n}\times\mathbb{S}^{m}\times\mathbb{S}^{s}.\\
\hline
1.&\text{\tt Compute the system of trilinear equations of \autoref{pres-sys}}.\\
2.&\text{\tt Add the polynomial equations }\|x\|^2=\|y\|^2=\|z\|^2=1.\\
3.&\text{\tt Compute a G\"obner basis for the resulting system, }\mathcal{I}.\\
4.&\text{\tt Find a basis }B\,\,\text{\tt of }\C[x_1,\ldots,x_{n},y_1,\ldots,y_m,z_1,\ldots,z_s]/\mathcal{I}.\\
5.&\text{\tt Compute the multiplication matrix of }\ell.\\
6.&\text{\tt Return the magnitude of the maximum real eigenvalue.}
\end{array}$$
}

Let's give an implementation of the algorithm in Maple.
The code computes the maximum of a trilinear form over $\mathbb{S}^{n-1}\times\mathbb{S}^{m-1}\times\mathbb{S}^{s-1}$.
The reader may change the values of $n$, $m$ and $s$ and the trilinear form, to get different examples.

{\scriptsize
\begin{verbatim}
> restart;with(Groebner):with(linalg):
> n:=2:m:=2:s:=2:
> L:=6*x[1]*y[1]*z[1]+3*x[2]*y[1]*z[1]-6*x[1]*y[2]*z[1]+16*x[2]*y[2]*z[1]-
  14*x[1]*y[1]*z[2]-15*x[2]*y[1]*z[2]-11*x[1]*y[2]*z[2]+8*x[2]*y[2]*z[2];
> #Step 1 and 2
> J:={add(x[i]^2,i=1..n)-1,add(y[j]^2,j=1..m)-1,add(z[k]^2,k=1..s)-1,
> seq(seq(x[i]*diff(L,x[j])-x[j]*diff(L,x[i]),j=1..i-1),i=1..n),
> seq(seq(y[i]*diff(L,y[j])-y[j]*diff(L,y[i]),j=1..i-1),i=1..m),
> seq(seq(z[i]*diff(L,z[j])-z[j]*diff(L,z[i]),j=1..i-1),i=1..s)}:
> #Step 3
> G:=Basis(J,'tord'):
> #Step 4
> ns,rv:=NormalSet(G, tord):
> #Step 5
> mulMat:=evalm(evalf(MultiplicationMatrix(L,ns,rv,G,tord))):
> #Step 6
> max(op(map(abs,map(Re,{eigenvalues(mulMat)}))));
\end{verbatim}
}

The following is a table that shows the time, in seconds, used to compute the maximum.
In the first column appears different values of $(n,m,s)$, in the second, the time used to compute the Steps 1 through 4
and in the third, the total time of the algorithm. We run a Maple 11 session on a 2.1GHz CPU, with 2GB of memory,\\
\begin{center}
\begin{tabular}{r|c|l}
$(n,m,s)$&Steps 1-4&Total time\\ \hline
$(2,2,2)$&0.03&0.33\\
$(2,2,3)$&0.05&0.79\\
$(2,2,4)$&0.09&0.99\\
$(2,2,5)$&0.14&1.20\\
$(2,3,3)$&0.31&7.13\\
$(2,3,4)$&0.89&30.03\\
$(3,3,3)$&5.06&397.28\\
\end{tabular}
\end{center}
Note that the computation of the multiplication matrix using
Gr\"obner basis, requires most of the time.
In \cite{MR1022943}, the authors proposed an algorithm to compute 
the multiplication matrix of $\ell$ directly; without Gr\"obner basis.

\section*{Acknowledgments.}
I would like to thanks Federico Holik, for very fruitful discussions and for presenting me this subject.
Also to the reviewer for his detailed and useful suggestions.
This work was supported by CONICET, Argentina.

\section*{References.}
\nocite{MR2754382} \nocite{MR2546345} \nocite{MR544868}\nocite{MR2746624}
\bibliographystyle{siam}
\bibliography{../../citas}

\def\cprime{$'$} \def\cprime{$'$}
\begin{thebibliography}{10}

\bibitem{MR0344384}
{\sc T.~M. Apostol}, {\em Mathematical analysis}, Addison-Wesley Publishing
  Co., Reading, Mass.-London-Don Mills, Ont., second~ed., 1974.

\bibitem{MR1022943}
{\sc W.~Auzinger and H.~J. Stetter}, {\em An elimination algorithm for the
  computation of all zeros of a system of multivariate polynomial equations},
  in Numerical mathematics, {S}ingapore 1988, vol.~86 of Internat.
  Schriftenreihe Numer. Math., Birkh\"auser, Basel, 1988, pp.~11--30.

\bibitem{BKP11}
{\sc G.~Ballard, T.~Kolda, and T.~Plantenga}, {\em Efficiently computing tensor
  eigenvalues on a gpu}, 2011 IEEE International Symposium on Parallel and
  Distributed Processing Workshops and PhD Forum,, 0 (2011), pp.~1340--1348.

\bibitem{MR0435072}
{\sc D.~N. Bernstein}, {\em The number of roots of a system of equations},
  Funkcional. Anal. i Prilo\v zen., 9 (1975), pp.~1--4.

\bibitem{Cartwright2011}
{\sc D.~Cartwright and B.~Sturmfels}, {\em The number of eigenvalues of a
  tensor}, Linear Algebra and its Applications,  (2011), pp.~--.

\bibitem{MR1931400}
{\sc P.~Comon}, {\em Tensor decompositions: state of the art and applications},
  in Mathematics in signal processing, {V} ({C}oventry, 2000), vol.~71 of Inst.
  Math. Appl. Conf. Ser. New Ser., Oxford Univ. Press, Oxford, 2002, pp.~1--24.

\bibitem{MR1780272}
{\sc L.~De~Lathauwer, B.~De~Moor, and J.~Vandewalle}, {\em A multilinear
  singular value decomposition}, SIAM J. Matrix Anal. Appl., 21 (2000),
  pp.~1253--1278 (electronic).

\bibitem{MR2447444}
{\sc V.~de~Silva and L.-H. Lim}, {\em Tensor rank and the ill-posedness of the
  best low-rank approximation problem}, SIAM J. Matrix Anal. Appl., 30 (2008),
  pp.~1084--1127.

\bibitem{MR2161984}
{\sc A.~Dickenstein and I.~Z. Emiris}, eds., {\em Solving polynomial
  equations}, vol.~14 of Algorithms and Computation in Mathematics,
  Springer-Verlag, Berlin, 2005.
\newblock Foundations, algorithms, and applications.

\bibitem{MR2648690}
{\sc T.-C. Dinh and N.~Sibony}, {\em Dynamics in several complex variables:
  endomorphisms of projective spaces and polynomial-like mappings}, in
  Holomorphic dynamical systems, vol.~1998 of Lecture Notes in Math., Springer,
  Berlin, 2010, pp.~165--294.

\bibitem{MR1644323}
{\sc W.~Fulton}, {\em Intersection theory}, vol.~2 of Ergebnisse der Mathematik
  und ihrer Grenzgebiete. 3. Folge. A Series of Modern Surveys in Mathematics
  [Results in Mathematics and Related Areas. 3rd Series. A Series of Modern
  Surveys in Mathematics], Springer-Verlag, Berlin, second~ed., 1998.

\bibitem{MR1264417}
{\sc I.~M. Gel{\cprime}fand, M.~M. Kapranov, and A.~V. Zelevinsky}, {\em
  Discriminants, resultants, and multidimensional determinants}, Mathematics:
  Theory \& Applications, Birkh\"auser Boston Inc., Boston, MA, 1994.

\bibitem{MR1182558}
{\sc J.~Harris}, {\em Algebraic geometry}, vol.~133 of Graduate Texts in
  Mathematics, Springer-Verlag, New York, 1992.
\newblock A first course.

\bibitem{MR0463157}
{\sc R.~Hartshorne}, {\em Algebraic geometry}, Springer-Verlag, New York, 1977.
\newblock Graduate Texts in Mathematics, No. 52.

\bibitem{MR2854605}
{\sc T.~G. Kolda and J.~R. Mayo}, {\em Shifted power method for computing
  tensor eigenpairs}, SIAM J. Matrix Anal. Appl., 32 (2011), pp.~1095--1124.

\bibitem{1574201}
{\sc L.-H. Lim}, {\em Singular values and eigenvalues of tensors: a variational
  approach}, in Computational Advances in Multi-Sensor Adaptive Processing,
  2005 1st IEEE International Workshop on, dec. 2005, pp.~129 --132.

\bibitem{MR2546345}
{\sc C.~Ling, J.~Nie, L.~Qi, and Y.~Ye}, {\em Biquadratic optimization over
  unit spheres and semidefinite programming relaxations}, SIAM J. Optim., 20
  (2009), pp.~1286--1310.

\bibitem{MR1720109}
{\sc A.~McLennan}, {\em The maximum number of real roots of a multihomogeneous
  system of polynomial equations}, Beitr\"age Algebra Geom., 40 (1999),
  pp.~343--350.

\bibitem{MR1342295}
{\sc A.~P. Morgan, A.~J. Sommese, and C.~W. Wampler}, {\em A
  product-decomposition bound for {B}ezout numbers}, SIAM J. Numer. Anal., 32
  (1995), pp.~1308--1325.

\bibitem{MR2296920}
{\sc G.~Ni, L.~Qi, F.~Wang, and Y.~Wang}, {\em The degree of the
  {E}-characteristic polynomial of an even order tensor}, J. Math. Anal. Appl.,
  329 (2007), pp.~1218--1229.

\bibitem{MR2270090}
{\sc L.~Qi}, {\em Eigenvalues and invariants of tensors}, J. Math. Anal. Appl.,
  325 (2007), pp.~1363--1377.

\bibitem{MR2746624}
{\sc B.~Savas and L.-H. Lim}, {\em Quasi-{N}ewton methods on {G}rassmannians
  and multilinear approximations of tensors}, SIAM J. Sci. Comput., 32 (2010),
  pp.~3352--3393.

\bibitem{MR544868}
{\sc L.~T. Watson}, {\em A globally convergent algorithm for computing fixed
  points of {$C^{2}$} maps}, Appl. Math. Comput., 5 (1979), pp.~297--311.

\bibitem{MR2754382}
{\sc X.~Zhang, C.~Ling, and L.~Qi}, {\em Semidefinite relaxation bounds for
  bi-quadratic optimization problems with quadratic constraints}, J. Global
  Optim., 49 (2011), pp.~293--311.

\end{thebibliography}
\end{document}